\documentclass [12pt]{article}
\usepackage{amsmath, amsfonts, amssymb, calrsfs}
\usepackage[russian]{babel}
\usepackage {graphicx}
\usepackage{wrapfig}
\usepackage{bbm}
\usepackage{mathrsfs}

  \DeclareMathAlphabet{\mathbbm}{U}{bbold}{m}{n}

\begin{document}

\makeatletter

\renewcommand{\thefootnote}{\fnsymbol{}}

\newcommand{\Author}[2]{{
\begin{center}\textbf{#1}\footnote{$\copyright$\ \Year\ #2}\end{center} \medskip}
    \renewcommand{\@evenhead}
    {\raisebox{0pt}[\headheight][0pt]%
{\vbox{\hbox to\textwidth{\issu--\thepage \hfill\strut {\sl
#1}}\hrule}}} }

\newcommand{\shorttitle}[1]{
\renewcommand{\@oddhead}{\raisebox{0pt}[\headheight][0pt]%
{\vbox{\hbox to\textwidth{\strut {\sl #1}
\hfill\issu--\thepage}\hrule}}} }

\headsep=3mm
\renewcommand{\section}[1]{\medskip\begin{center}\textbf{#1}\end{center}}
\newcommand{\subsec}[1]{\par\smallskip{\bf #1}}
\newcommand{\subsubsec}[1]{\qquad\qquad{\bf #1}}
 \def\beginproof{\smallskip\text{$\vartriangleleft$}}
 \def\endproof{\text{$\vartriangleright$}}
 \def\Endproof{\text{$\vartriangleright$}\smallskip}

\makeatother

\newcommand{\beginarticle}{\newpage\thispagestyle{empty}
\begin{flushright}
\sl T« ¤ЁЄ ўЄ §бЄЁ© ¬ вҐ¬ вЁзҐбЄЁ© ¦га­ «\\
\moons, \Year, T®¬~\Vol, TvЇгбЄ~\issu
\end{flushright}
\bigskip}

\newcommand{\UDC}[1]{{\sl L-i \bf#1}\par\vspace{5mm}}
\newcommand{\Title}[1]{\begin{center}\uppercase{#1}\end{center} \vspace{0mm}}
\newcommand{\TITLE}[1]{\begin{center}{#1}\end{center} \vspace{0mm}}
\renewcommand{\title}[1]{\begin{center}#1\end{center} \vspace{0mm}}
\newcommand{\Abstract}[1]{\hangindent20pt\hangafter=0\noindent{\footnotesize
#1}\bigskip\par\medskip }

\newcommand{\bib}[2]{{\baselineskip=11pt\footnotesize\item{}
\textsl{#1\/}\ #2}}
\newcommand{\Address}[1]{\par\bigskip\baselineskip=11pt\hangindent24pt
\hangafter=0\noindent{\footnotesize#1}\par\vfill\eject \par \normalsize}
\newcommand{\Endarticle}[2]{\bigskip\bigskip{\textit{#1\hfill Tв вмп Ї®бвгЇЁ« 
#2}} \vfill\eject}

\newcommand{\Endproc}{\rm} 
\newcommand{\Proclaim}[1]{\smallskip{\textbf{#1\/~}}\sl}
\newcommand{\proclaim}[1]{{\textbf{#1\/~}}\sl}

\newcommand{\Teorema}[1]{\smallskip{\textbf{TҐ®аҐ¬ #1.\/~}}\sl}
\newcommand{\teorema}[1]{{\textbf{TҐ®аҐ¬ #1.\/~}}\sl}
\newcommand{\Sledstvie}[1]{\smallskip{\textbf{T«Ґ¤бвўЁҐ#1.\/~}}\sl}
\newcommand{\sledstvie}[1]{{\textbf{T«Ґ¤бвўЁҐ#1.\/~}}\sl}
\newcommand{\Lema}[1]{\smallskip\textbf{{TҐ¬¬ #1}.\/~}\sl}
\newcommand{\lema}[1]{\textbf{{TҐ¬¬ #1}.\/~}\sl}
\newcommand{\Predl}[1]{\smallskip\textbf{{iаҐ¤«®¦Ґ­ЁҐ#1}.\/~}\sl}
\newcommand{\predl}[1]{\textbf{{iаҐ¤«®¦Ґ­ЁҐ#1}.\/~}\sl}
\newcommand{\Utv}[1]{\smallskip\textbf{{LвўҐа¦¤Ґ­ЁҐ#1}.\/~}\sl}
\newcommand{\utv}[1]{\textbf{{LвўҐа¦¤Ґ­ЁҐ#1}.\/~}\sl}
\newcommand{\Def}[1]{\smallskip{\bf +ЇаҐ¤Ґ«Ґ­ЁҐ#1.}}
\newcommand{\Zam}[1]{\smallskip{\sc i ¬Ґз ­ЁҐ#1.}}
\newcommand{\zam}[1]{{\sc i ¬Ґз ­ЁҐ#1.}}
\newcommand{\Primer}[1]{\smallskip{\sc iаЁ¬Ґа#1.}}
\newcommand{\primer}[1]{{\sc iаЁ¬Ґа#1.}}
\newcommand{\Zadacha}[1]{\smallskip{\sc i ¤ з #1.}}
\newcommand{\zadacha}[1]{{\sc i ¤ з #1.}}
%
\newcommand{\Lemma}[1]{\smallskip\textbf{{Lemma#1}.\/~}\sl}
\newcommand{\lemma}[1]{\textbf{{Lemma#1}.\/~}\sl}
\newcommand{\Theorem}[1]{\smallskip\textbf{{Theorem#1}.\/ }\sl}
\newcommand{\theorem}[1]{{\textbf{Theorem#1.\/~}}\sl}
\newcommand{\Corollary}[1]{\smallskip\textbf{{Corollary#1}.\/~}\sl}
\newcommand{\Proposition}[1]{\smallskip\textbf{{Proposition#1}.\/~}\sl}
\newcommand{\defin}[1]{{\sc Definition#1.}}
\newcommand{\Definition}[1]{\smallskip{\sc Definition#1.}}
\newcommand{\Remark}[1]{\smallskip{\sc Remark#1.}~}

\renewcommand{\Re}{{\rm Re\,}}
\newcommand{\Div}{\text{div\,}}
\newcommand{\pd}[2]{\frac{\partial #1}{\partial #2}}
\newcommand{\pg}[2]{\frac{{\partial}^2 #1}{\partial #2^2}}
\newcommand{\const}{\textrm{const}}
\newcommand{\grad}{\textrm{grad\/}}
\newcommand{\sgn}{\mathop{\fam0 sgn}}
\newcommand{\rot}{\mathrm{rot}}
\renewcommand{\leq}{\leqslant}
\renewcommand{\geq}{\geqslant}
\renewcommand{\le}{\leqslant}
\renewcommand{\ge}{\geqslant}
\newcommand{\im}{\mathop{\fam0 Im}}
\newcommand{\cl}{\mathop{\fam0 cl}}
\newcommand{\mix}{\mathop{\fam0 mix}}
\newcommand{\ind}{\mathop{\fam0 ind}}
\newcommand{\proj}{\mathop{\fam0 proj}}

\def\reduce{\mskip-5mu }
\def\Bnorml{\mathopen{\kern1pt \vrule height8.2pt depth2.8pt width1pt\kern1.5pt}}
\def\Bnormr{\mathclose{\kern1.5pt \vrule height8.2pt depth2.8pt width1pt\kern1pt}}
\def\widevert{\kern1.5pt \vrule height7.5pt depth2.5pt width1pt\kern1pt}

\def\sfatop{\mathopen{\kern1pt \vrule height5.5pt depth1.5pt width1pt\kern1pt}}
\def\sfatcl{\mathclose{\kern1pt\vrule height5.5pt depth1.5pt width1pt\kern1pt}}

\def\leftbbr{\mathopen{\kern1pt \vrule height14pt depth8pt width1pt\kern1.5pt}}
\def\leftbbrr{\mathopen{\kern1pt\vrule height16pt depth13pt width1pt\kern1.5pt}}
\def\rightbbrr{\mathclose{\kern1.5pt\vrule height16pt depth13pt width1pt\kern1pt}}
\def\rightbbr{\mathclose{\kern1.5pt\vrule height14pt depth8pt width1pt\kern1pt}}


\newcommand{\shortpage}{\enlargethispage{-\baselineskip}}
\newcommand{\largepage}{\enlargethispage{+\baselineskip}}




 \def\Orth{\mathop{\fam0 Orth}}
 \def\Otimes{\overset{\underline{\phantom{\ \,}}}\otimes}
 \def\cone{\mathop{\fam0 cone}}
 \def\co{\mathop{\fam0 co}}
 \def\fin{\mathop{\fam0 fin}}
 \def\im{\mathop{\fam0 im}}
 \def\Prt{\mathop{\fam0 Prt}}
 \def\End{\mathop{\fam0 End}}
 \setcounter{equation}{0}

 \def\[{\mathopen{\kern1pt\/ \vrule height9pt depth2.5pt width1pt\kern1.5pt}}
 \def\]{\mathclose{\kern1.5pt\vrule height9pt depth2.5pt width1pt\kern1pt}}
 \def\sfatop{\mathopen{\kern1pt \vrule height5.5pt depth1.5pt width1pt\kern1pt}}
 \def\sfatcl{\mathclose{\kern1pt\vrule height5.5pt depth1.5pt width1pt\kern1pt}}

 \def\tv@rt{{\vert\mkern-2mu\vert\mkern-2mu\vert}}
 \def\tvert#1\tvert{\mathopen\tv@rt#1\mathclose\tv@rt}
 \def\beginproof{\par\mbox{$\vartriangleleft$}}
 \def\endproof{\text{$\vartriangleright$}}


 \let\thefootnote\relax

 \Title{Geometric characterization of \\preduals of injective Banach lattices}


 \begin{center}
 Kusraev A.~G. and Kutateladze~S.~S.
 \end{center}

 {\bf Abstract}{ The paper deals with the study of Banach spaces whose duals are injective Banach lattices.  Davies in 1967 proved
 that an ordered Banach space is an $L^1$-predual space if and only if it is a simplex space. In 2007 Duan and Lin
 proved that a real Banach space is an $L^1$-predual space if and only if its every four-point subset is centerable.
 We prove the counterparts of these remarkable results for injectives by the new machinery of Boolean
 valued transfer from $L^1$-spaces to injective Banach lattices.\par}

 {\hangindent17pt\hangafter=0\noindent\footnotesize
 \par}

 \setcounter{equation}{0}

 \section{1. Introduction}

  A real Banach lattice $X$ is said to be \textit{injective\/} if, given a Banach lattice $Y$, a closed vector sublattice
  $Y_0\subset Y$, and a positive linear operator $T_0:Y_0\to X$, there exists a positive linear extension
 $T:Y\to X$ with $\|T_0\|=\|T\|$. Equivalently, $X$ is an injective Banach lattice if, whenever $X$ is lattice
 isometrically embedded into a Banach lattice $Y$, there exists a positive contractive projection from $Y$ onto
 $X$. Thus, the injective Banach lattices are the injective objects in the category of Banach lattices with positive
 contractions as morphisms. The first example of an injective Banach lattice was given by Abramovich \cite{Abr}, while a
 systematic study was started by Lotz \cite{Lotz}, who introduced the term ``injective Banach lattice.''
 A remarkable contribution
 to the study of injective Banach lattices was made by Cartwright \cite{Car} and Haydon \cite{Hay} who discovered important
 geometric and analytical properties of injective Banach lattices.

 In  \cite[Theorems 4.1 and 4.4]{Kus_IBL}, the \textit{transfer principle from $AL$-spaces to injective Banach spaces\/} was
 established which can be formulated as follows: (1)~Every injective Banach lattice embeds into an appropriate Boolean
 valued model of set theory, becoming an $AL$-space, (2)~Each theorem about $AL$-spaces within ${\rm ZFC}$ ($=$~Zermelo--Fraenkel set
 theory with the axiom of choice) has a counterpart for the original injective Banach lattice interpreted as a Boolean valued
 $AL$-space; (3)~Translation of theorems from $AL$-spaces to injective Banach lattices is carried out by appropriate general
 operations and principles of Boolean valued analysis.


 One of the intriguing problems that dated from the works of Grothendieck~\cite{Gro} and Lindenstrauss \cite{Lin} is to~describe the Banach spaces whose duals are isometric (isomorphic) to an $AL$-spaces. The injective version
 of this problem was posed in \cite[Problem 5.16]{KW}:~\textit{Classify and characterize the Banach spaces whose duals are
 injective Banach lattices}. The present article can be considered as an attempt to pave a way towards the study of the problem.

 We organize the article as follows: Section 2 collects the needed information about Boolean algebras of projections in
 real Banach spaces. The key fact in our study, which is due to Cunningham \cite{Cun1}, claims that the collection of
 $L$-projections in a Banach space is a Bad\'e complete Boolean algebra (Theorem 2.10). In Section 3, we outline the
 approach of Boolean valued analysis to the characterization of the Banach spaces predual to injective Banach lattices.
 We demonstrate that a Banach space admits a Boolean valued representation which is a Banach space without
 nontrivial $L$-projections, while the dual space is represented as the dual Banach space without nontrivial
 $M$-projections (Theorem 3.8). In Section 4 we present some conditions on the  unit ball of an ordered Banach
 space that are necessary and sufficient for the space to be a predual to an injective Banach lattice (Theorem 4.9).
 This is carried out by Boolean valued interpretation of the result by Davies \cite{Dav} which asserts
 that an ordered Banach space is $X$ an $L^1$-predual space if and only if $X$ is a simplex space. In Section 5
 we characterize the Banach spaces preduals to injective Banach lattices in terms of centerable sets (Theorem 5.6).
 To this end we interpret, in a Boolean valued model, the following result by Duan and Lin \cite{DL}: A Banach space $X$
 is $L^1$-predual if and only if every four-point  subset of $X$ is centerable.

 For the theory of Banach lattices and positive operators we refer to the books \cite{LuZ, MN}. The needed
 information on the theory of Boo\-lean valued models of set theory is briefly presented in \cite[Chapter 9]{DOP} and
 \cite[Chapter 1]{KKTop}; details can be found in~\cite{Bell, BVA}.

 Throughout the sequel $\mathbb{B}$ is a complete Boolean algebra with join $\vee$, meet $\wedge$, complement
 $(\cdot)^\ast$, unit (top) $\mathbbm{1}$, and zero (bottom) $\mathbb{O}$; while $\mathbb{P}(E)$ stands for the Boolean
 algebra of band projections in a vector lattice $E$.
 A~{\it partition of unity\/} in~$\mathbb{B}$ is a~family $(b_{\xi})_{\xi\in \Xi}\subset\mathbb{B}$ such that
 $\bigvee_{\xi\in \Xi} b_{\xi}=\mathbbm{1}$ and $b_{\xi}\wedge b_{\eta}=\mathbb{O}$ whenever $\xi\ne \eta$. We let $:=$
 denote the assignment by definition, while $\mathbb{N}$ and $\mathbb{R}$ symbolize the naturals and the reals.

 \section{2. Preliminaries}

 In what follows, we need some information about Boolean algebras of projections
  in real vector spaces. When speaking
 of a {\it Boolean algebra of projections\/} in a~vector space $X$ we always mean some~set $\mathcal{P}$ of commuting
 idempotent linear operators; i.e. projections, on $X$ which is a Boolean algebra under the operations
 $$
 \pi \wedge \rho\!:=\pi \circ \rho=\rho \circ \pi, \quad \pi \vee
 \rho=\pi + \rho - \pi \circ \rho,      \quad \pi^\ast=I_X-\pi \quad
 (\pi,\rho \in\mathcal{P})
 $$
 and in which the~zero and the identity operators in $X$ serve as the~top and bottom of $\mathcal P$.

 \Definition{~2.1}~Let $\mathcal{P}_X$ and $\mathcal{P}_Y$ be Boolean algebras of projections in $X$ and $Y$, respectively,
 both isomorphic to $\mathbb{B}$. An operator $T:X\rightarrow Y$ is called $\mathbb{B}$-{\it linear}, if $T$ is linear and
 $\psi(b)\circ T=T\circ \varphi(b)$ for all $b\in\mathbb{B}$, where $\varphi$ and $\psi$ are Boolean
 isomorphisms of
 $\mathbb{B}$ onto $\mathcal{P}_X$ and $\mathcal{P}_Y$, respectively. A one-to-one $\mathbb{B}$-linear operator is called
 a~$\mathbb{B}$-{\it isomorphism}, and an isometric $\mathbb{B}$-{\it isomorphism} of Banach spaces is called
 a~{\it $\mathbb B$-isometry.} A~$\mathbb B$-isometric lattice homomorphism between Banach lattices
 is referred to as a~\textit{lattice $\mathbb{B}$-isometry}.
 \smallskip

 In the sequel, we identify $\mathbb{B}$ with $\mathcal{P}_X$ and $\mathcal{P}_Y$ and so we can
 say that $T$ commutes with all projections from $\mathbb{B}$ and write $b\circ T=T\circ b$ for $b\in\mathbb{B}$.
 Let us consider some situations in which a Boolean algebra of projections is associated with a norm. To start, take
 a vector space $X$ and a vector lattice $E$.

 \Definition{~2.2}~An~{\it $E$-valued norm\/} is a~mapping $\[\cdot\]:X\to E_+$ such that~$\[x\] =0$ implies that $x=0$,
 while $\[\lambda x\] =|\lambda|\[x\] $ and $\[x+y\]\leq \[ x\] +\[ y\] $ for all $x,y\in X$ and $\lambda\in\mathbbm{R}$.
 A \textit{lattice normed space\/} over $E$ is a pair $(X,\[\cdot\])$, where $\[\cdot\]$ is an $E$-valued norm on
 $X$. An~$E$-valued norm (as well as $X$ itself) is said to be {\it decomposable\/} if, for each decomposition
 $\[ x\] =e_1+e_2$ with $e_1,e_2\in E_+$ and $x \in X$, there exist $x_1,x_2 \in X$ such that $x=x_1+x_2$ and
 $\[ x_k\] =e_k$ $(k \!:= 1,2)$. Put $\[A\]\!:=\{\[x\]:\,x\in A\}$.

 \Proposition{~2.3 {\rm \cite[2.1.3]{DOP}}}Assume that $X$ is a decomposable lattice normed space over a vector lattice $E$ with the
 projection property and $E=\[X\]^{\perp\perp}$. Then there is a~complete Boolean algebra $\mathcal{P}(X)$ of
 projections in $X$ and a~Boolean isomorphism $h$ from $\mathbb{P}(E)$ onto $\mathcal{P}(X)$ such that
 $b\[x\]=\[h(b)x\]$ for all $b\in\mathbb{P}(E)$ and $x\in X$. Moreover, for $\pi\!:=h(b)$ we have
 \begin{equation}\label{D}
 \[\pi x+\pi^\ast y\]=\pi\[x\]+\pi^\ast\[y\]\quad(x,y\in X).
 \end{equation}
 \Endproc

 \Definition{~2.4}~A sequence $(x_n)$ in $X$ is said to be $E$-\textit{uniformly convergent} to $x\in X$ (respectively,
 $E$-\textit{uniformly Cauchy}) whenever the sequence $(\[x-x_n\])$ is uniformly convergent to zero (respectively,
 $(\[x_n-x_m\])$ is uniformly Cauchy) in $E$.  Say that $X$ is $E$-\textit{uniformly complete}  whenever every
 $E$-uniformly Cauchy sequence in $X$ is $E$-uniformly convergent.

 \Definition{~2.5}~Consider a lattice normed space $(X,\[\cdot\])$ over a Banach lattice $E$ and put
 $\mathbb{B}\!:=\mathbb{P}(E)$. Endow $X$ with the \textit{mixed norm} $\|\cdot\|$ defined as $\|x\|\!:=\big\|\[x\]\big\|_E$
 $(x\in X)$. A subset $S\subset X$ is called $\mathbb{B}$-\textit{bounded\/} if there exists $a\in X$ such that
 $\|\pi x\|\leq\|\pi a\|$ for all $x\in S$ and $\pi\in\mathbb{B}$. Say that $X$ is $\mathbb{B}$-\textit{complete\/} whenever,
 given a~partition of unity $(b_\xi)$ in $\mathbb{B}$ and a $\mathbb{B}$-bounded family $(x_\xi)$ in $X$, there exists a
 unique $x\in X$ such that $b_\xi x=b_\xi x_\xi$ for all $\xi\in\Xi$.

 \Proposition{~2.6 {\rm \cite[7.2.2]{DOP}}}~Let~$(X,\[\cdot\])$ be a lattice normed space over a~Banach lattice $E$ and let
 $\|\cdot\|$ stand for the associated mixed norm. Then $(X,\|\cdot\|)$ is a~Banach space if and only if $(X,\[\cdot\])$ is
 $E$-uniformly complete.
 \Endproc

 \Definition{~2.7}~A~{\it Ba\-nach--Kantorovich space\/} over a vector lattice $E$ is a vector space $X$ with a
 decomposable norm $\[\cdot \]:X\rightarrow E$ which is norm complete in the sense that, given a net
 $(x_\alpha)_{\alpha\in\mathrm{A}}$ in $X$ with
 $(\[ x_\alpha-x_\beta\] )_{(\alpha,\beta)\in \mathrm{A}\times\mathrm{A}}$ order convergent to the  zero
 of $E$, there exists $x\in X$ such that $(\[ x_\alpha-x\])_{\alpha\in\mathrm{A}}$ is order convergent
 to the zero of $E$.

 \Proposition{~2.8}A decomposable lattice normed space $(X,\[\cdot\])$ over a~Banach lattice~$E$ is a
 Banach--Kantorovich space if and only if the associated mixed norm space $(X,\|\cdot\|)$ is a~$\mathbb{B}$-complete
 Banach space.
 \Endproc

 \beginproof~Observe first that the condition $\|\pi x\|\leq\|\pi a\|$ for all $\pi\in\mathbb{B}$ is equivalent to the
 inequality $\[x\]\leq\[a\]$. Indeed, $\[x\]\leq\[a\]$ obviously implies that
 $\|x\|=\big\|\[a\]\big\|_E\leq\big\|\[x\]\big\|_E=\|a\|$ by monotonicity of $\|\cdot\|_E$. If $\[x\]\leq\[a\]$
 is false then, using Proposition 2.3, we can pick up the projection $\pi_0\in\mathbb{B}$ and a number $\varepsilon>0$
 so that $\[\pi_0 x\]\geq(1+\varepsilon)\[\pi_0 a\]$. It follows that $\|\pi_0 x\|\geq(1+\varepsilon)\|\pi_0 a\|>\|\pi_0 a\|$;
 a contradiction. Now it is clear that a set $A\subset X$ is $\mathbb{B}$-bounded if and only if $\[A\]$ is order bounded
 in $E$. It remains to recall that a decomposable lattice normed space over $E$ is norm complete if and only if it is
 disjointly complete and uniformly $E$-complete (see \cite[5.3.5 and 5.4.7]{BVA}).~\endproof

 We now consider special projections in a Banach space whose presence indicates that the Banach space has some
 features of an $AL$-space or $AM$-space.


 \Definition{~2.9}~A~projection $\pi$ on a Banach space $X$ is said to be an $M$-\textit{projection} if
 $\|x\|=\|\pi x\|+\|x-\pi x\|$ for all $x\in X$ and an $L$-\textit{projection} if $\|x\|=\max\{\|\pi x\|, \|x-\pi x\|\}$
 for all $x\in X$. The $L$-projections and $M$-projections different from the zero and the identity are referred to as
 \textit{nontrivial}. The sets of all $L$-projections and $M$-projections on $X$ will be  denoted by $\mathbb{P}_L(X)$
 and $\mathbb{P}_M(X)$, respectively.


 A simple induction argument shows that for all $x\in X$ every finite collection of pairwise disjoint projections
 $\{\pi_1,\dots,\pi_n\}$ in $\mathbb{P}_L(X)$ (respectively, $\mathbb{P}_M(X)$) with $\pi_0=\pi_1+\dots+\pi_n$ we have
 \begin{equation}\label{LM}
 \|\pi_0x\|\!=\|\pi_1x\|+\dots+\|\pi_nx\|\quad(\text{respectively,~~}\|\pi_0x\|\!=\max\{\|\pi_1x\|,\dots,\|\pi_nx\|\}).
 \end{equation}

 The following result is due to Cunningham \cite[Theorem 2.5]{Cun1} and \cite[Theorem 5]{Cun2}.

 \Theorem{~2.10}~For a Banach space $X$ the following hold:

 $(1)$~$\mathbb{P}_L(X)$ is a complete Boolean algebra.

 $(2)$~$\mathbb{P}_M(X)$ is a $($generally not complete$)$ Boolean algebra.

 $(3)$~If $X'$ is the dual of $X$ then $\mathbb{P}_M(X')$ is isomorphic to $\mathbb{P}_L(X)$. 
 \Endproc

 \beginproof~The proof can be found in \cite[Theorem 1.10]{HWW}. \endproof

 It follows that the Boolean algebra $\mathbb{P}_M(X')$ is complete. Another examples of complete Boolean algebras of
 $M$-projections provide injective Banach lattices: The Boolean algebra $\mathbb{P}_M(X)$ is complete for an arbitrary
 injective Banach lattice $X$; moreover the unit ball of $X$ is $\mathbb{P}_M(X)$-complete; see Lemma 4.5 below.

 Cunningham \cite[Lemma 2.4]{Cun1} proved also that $\mathbb{P}_L(X)$ is \textit{Bad\'e complete:}
 If $(\pi_\alpha)$ is an increasing net of $L$-projections and $\pi\!:=\sup_\alpha\pi_\alpha$
 then $(\pi_\alpha x)$ is norm convergent to $\pi x$ for all $x\in X$. It follows that for every
 $0\ne x\in X$ the function $\mu_x:\mathbb{P}_L(X)\to\mathbb{R}$ defined as $\mu_x:\pi\mapsto\|\pi x\|$
 is a nonzero order continuous measure. Hence, $\mathbb{P}_L(X)$ has a separating set of order continuous measures
 $($normal measure$)$ or, equivalently, the  representation space of $\mathbb{P}_L(X)$ is \textit{hyperstone}.


 \Lemma{~2.11}~Let $X$ be a Banach space. Assume that $a\in X$ and  $(x_\xi)_{\xi\in\Xi}$ in $X$ satisfies the
 condition $\|\pi x_\xi\|\leq\|\pi a\|$ for all $\xi\in\Xi$ and $\pi\in\mathbb{P}_L(X)$. Then for every partition of unity
 $(\pi_\xi)_{\xi\in\Xi}$ in $\mathbb{P}_L(X)$ there exists a unique $x\in X$ such that $\pi_\xi x_\xi=\pi_\xi x$ for all
 $\xi\in\Xi$. Moreover, $x=\sum_{\xi\in\Xi}\pi_\xi x_\xi$ and $\|\pi x\|\leq\|\pi a\|$ for all $\pi\in\mathbb{P}_L(X)$.
 \Endproc

 \beginproof~Suppose that $(x_\xi)$ and $a\in X$ meet the hypotheses of the lemma. Let $\Theta$ be the set of all finite
 subsets of $\Xi$. Given $\theta\in\Theta$, put
 $$
 y_\theta\!:=\sum_{\xi\in\theta}\pi_\xi x_\xi,\quad
 \sigma_\theta\!:=\bigvee_{\xi\in\theta}\pi_\xi,\quad
 \sigma\!:=\bigvee_{\xi\in\Xi}\pi_xi,\quad
 \rho_\theta  \!:=\sigma -\sigma_\theta.
 $$
 Take $\theta,\theta_1,\theta_2\in\Theta$ with $\theta\subset\theta_1\cap\theta_2$ and denote by $\theta^\prime$ and
 $\theta_1\triangle\theta_2$ the complement of $\theta$ and the symmetric difference of $\theta_1$ and $\theta_2$,
 respectively. Now using \eqref{LM}, we have
 $$
 \|y_{\theta_1}-y_{\theta_2}\|
 =\bigg\|\sum_{\xi\in\theta_1\triangle\theta_2}\pi_\xi x_\xi\bigg\|
 =\sum_{e\in\theta_1\triangle\theta_2}\|\pi_\xi x_\xi\|\leq\sum_{e\in\theta^\prime}\|\pi_\xi a\|
 =\|\rho_\theta a\|.
 $$
 By hypothesis $(\delta _\theta)_{\theta\in\Theta}$ decreases to zero, so that $(y_\theta)_{\theta\in\Theta}$ is Cauchy and
 there exists $x\!:=\lim_{\theta\in\Theta}y_\theta=\sum_{\xi\in\Xi}\pi_\xi x_\xi$ in $X$. If $\xi\in\theta$ then
 evidently $\pi_\xi y_\theta=\pi_\xi x_\xi$ and passage to the limit yields $\pi_\xi x=\pi_\xi x_\xi$. Moreover,
 $\|\pi\pi_\xi y_\theta\|\leq \|\pi a\|$ by hypothesis and so $\|\pi\sigma_\theta x\|\leq \|\pi a\|$ for all
 $\theta\in\Theta$; again, passage to the limit in the inequality yields  $\|\pi x\|\leq\|\pi a\|$.~\endproof

 \Remark{~2.12}$L$-projections and $M$-projections were first studied by Cunninghem \cite{Cun1, Cun2}.
 He demonstrated in \cite{Cun2} that a Banach space can be represented as a section space of the~Banach bundle
 whose fibres have no nontrivial $L$-projections. According to Theorem 2.10 $L$-projections and $M$-projections are
 mutually dual. Nevertheless, there is a striking asymmetry in these dualities, because the collection
  of
 $M$-projections may be too scarce; for example, there exist no nontrivial $M$-projections on $C([0,1])$. A natural
 generalization of the concepts of $L$-projection and $M$-projection is the concept of $L^p$-projection; see \cite{Beh_ETC, HWW}.

 \section{3. Boolean Valued Representation}

 In this section we outline the approach of Boolean valued analysis  to the problem mentioned in the Introduction.
 Boolean valued analysis is the technique of studying properties of an arbitrary mathematical object by comparison
 between its representations in two different Boolean valued models of set theory. As the models, we usually take
 the \textit{von Neumann universe} ${\mathbb{V}}$ (the mundane embodiment of the classical Cantorian paradise) and
 the \textit{Boolean valued universe} $\mathbb{V}^{({\mathbb B})}$ (a specially-trimmed universe whose construction
 utilizes a complete Boolean algebra $\mathbb{B}$).  The pair $(\mathbb{V},\in)$ is a standard model of $\rm ZFC$. The
 principal difference between $\mathbb{V}$ and
 $\mathbb{V}^{({\mathbb B})}$ is the way of verification of statements. In fact, there is a natural way of assigning to each
 statement $\phi$ about $x_1,\dots,x_n\in\mathbb{V}^{({\mathbb B})}$ the \textit{Boolean truth-value\/}
 $[\![\phi(x_1,\dots,x_n)]\!]\in\mathbb{B}$. The sentence $\phi(x_1,\dots,x_n)$ is called
 {\it true within
 $\mathbb{V}^{(\mathbb{B})}$} if $[\![\phi(x_1,\dots,x_n)]\!]=\mathbbm{1}$. For every complete Boolean algebra $\mathbb{B}$,
 all the theorems of $\rm ZFC$ are true in $\mathbb{V}^{(\mathbb{B})}$. There is a~smooth mathematical
 technique for interplay between the interpretations of one and the same fact in the two models~${\mathbb{V}}$
 and~${\mathbb{V}}^{({\mathbb B})}$. The relevant {\it ascending-and-descending machinery\/} rests on the functors of
 \textit{canonical embedding} $X\mapsto X^{\scriptscriptstyle\wedge}$ and  \textit{ascent} $X\mapsto X{\uparrow}$, both
 acting from $\mathbb{V}$ to $\mathbb{V}^{({\mathbb B})}$, and the functor of \textit{descent} $X\mapsto X{\downarrow}$, acting from
 $\mathbb{V}^{({\mathbb B})}$ to $\mathbb{V}$; see \cite{BVA, KKTop} for details.
 \smallskip

 These functors are applicable, in particular, to algebraic structures. Applying the Transfer and Maximum Principles
 to the $\rm ZFC$-theorem on the existence of the reals, find $\mathcal{R}\in \mathbb{V}^{(\mathbb{B})}$, called the
 \textit{reals within} $\mathbb{V}^{(\mathbb{B})}$ satisfying $[\![\mathcal{R}\text{~~is the reals}]\!]= \mathbbm{1}$
 and $[\![1^{\scriptscriptstyle\wedge}\in\mathbb{R}^{\scriptscriptstyle\wedge}\subset\mathcal{R}]\!]= \mathbbm{1}$.
 The \textit{Gordon Theorem\/} \cite{Gor1} states that the descent $\mathcal{R}{\downarrow}$ of $\mathcal{R}$ (with
 the descended operations and order) is a universally complete vector lattice. The mapping
 $\chi:\mathbb{B}\to \mathbb{P}(\mathrsfs R{\downarrow})$ is defined by putting
 $\chi (b) x\!:=\operatorname{mix} \{bx,b^*0\}$ for $x\in \mathrsfs R{\downarrow}$ and $b\in \mathbb{B}$. In more detail,
 $\chi$ is uniquely determined by the  relations%
 \begin{equation}\label{CHI}
 \chi(b)\in\mathbb{P}(\mathcal{R}{\downarrow}),\quad
 b\le [\![\,\chi (b) x=x\,]\!],\quad
 b^*\le [\![\,\chi (b)x=0\,]\!].
 \end{equation}
 Then $\chi$ is a Boolean isomorphism from
 $\mathbb{B}$ onto $\mathbb{P}(\mathcal{R}{\downarrow})$ such that
 \begin{equation}\label{BV}
 \chi (b) x=\chi (b) y\Longleftrightarrow b\leq[\![\,x=y\,]\!],\quad
 \chi (b) x\le \chi (b) y\Longleftrightarrow b\le [\![\,x\leq y\,]\!]
 \end{equation}
 for all $x,y\in\mathcal{R}$ and $b\in\mathbb{B}$; see \cite[Theorem 5.2.2]{BVA} and \cite[Theorem 2.2.4]{KKTop}. Moreover,
 the universally complete vector lattice~$\mathcal{R}{\downarrow}$ endowed with the descended multiplications is a~semiprime
 $f$-algebra with the order and ring unit~$\mathbf{1}\!:=1^{\scriptscriptstyle\wedge}$. Moreover, for every~$b\in\mathbb{B}$
 the band projection $\chi(b)$ acts as multiplication by the $\chi (b)\mathbf{1}$; see \cite[Theorem 2.3.2]{KKTop}.

 \Definition{~3.1}~Consider a~Banach space $(\mathcal X,\rho)$ in $\mathbb{V}^{(\mathbb{B})}$ and an order dense ideal $E$ in
 $\mathcal{R}{\downarrow}$. The descent $\[{\cdot}\]\!:=\rho{\downarrow}({\cdot})$ of $\rho:\mathcal{X}\to\mathcal{R}$ is
 a mapping from $\mathcal{X}{\downarrow}$ to $\mathcal{R}{\downarrow}$ and we can define
 the $E$-\textit{descent} of $\mathcal{X}$ as the subspace
 $\mathcal{X}{\downarrow}\!^E\!:=\{x\in\mathcal{X}{\downarrow}:\, \[x\]\in E\}$ of $\mathcal{X}{\downarrow}$ with the induced $E$-valued norm.
  \smallskip

 {\bf Theorem~3.2.}~{\sl If $\mathcal{X}$ is a~Banach space within $\mathbb{V}^{(\mathbb{B})}$, then the $E$-descent
 $X\!:=\mathcal{X}{\downarrow}\!^E$ is a Banach--Kantorovich space over $E$
 and $X$ can be endowed with the structure of a faithful unital module over the $f$-algebra $\Orth(E)$ so that
 $b\leq [\![\,x=0\,]\!]$ if and only if $\chi(b)x=0$ for all
 $b\in \mathbb{B}$ and $x\in X$. Conversely, if $X$ is a Banach--Kantorovich space over $E$ with $E=\[X\]^{\perp\perp}$
 and\/ $\mathbb{B}\!:=\mathbb{P}(E)$, then there exists a Banach space $\mathcal{X}$ unique up to linear isometry
 within $\mathbb{V}^{(\mathbb{B})}$ whose $E$-descent is linearly isometric to $X$ $($in the sense of $E$-valued norms$)$.}

 \beginproof~The proof can be extracted from  \cite[Theorems 5.4.1 and 5.4.2]{BVA}.~\endproof

 It is worth to note the two particular cases of Theorem 3.2 which concern $AL$-spaces and $AM$-spaces.
 Let $M$ be the {\it bounded part\/} of the universally complete vector lattice $\mathcal{R}{\downarrow}$, i.e. $M$
 consists of all $x\in\mathcal{R}{\downarrow}$ with $|x|\leq C\mathbf{1}$ for some $C\in\mathbb{R}$. Endow $M$ with the
 norm $\|m\|_\infty\!:=\inf\{0<\lambda\in{\mathbb R}:\,|m|\leq\lambda\mathbf{1}\}$. Then $M$ is a Dedekind complete
 $AM$-space with unit and $\mathbb{B}=\mathbb{P}(M)$. Putting $E\!:=M$ in Theorem 3.2 we arrive at the following:
 The $M$-descent $X\!:=\mathcal{X}{\downarrow}\!^M$ endowed with the mixed norm $\|x\|\!:=\big\|\[x\]\big\|_\infty$ is
 a Banach space; see Proposition 2.8. The $M$-descent is also called the \textit{bounded descent} and still denoted by
 $\mathcal{X}{\Downarrow}$; see \cite{BVA, KKTop}.

 \Definition{~3.3}~The $\mathbb{B}$-dual $X^{\scriptscriptstyle\#}$ of a lattice normed space $X$ over $E$ is defined
 as the lattice normed space over $\Orth(E)$, where $X^{\scriptscriptstyle\#}$ consists of all linear operators
 $x^{\scriptscriptstyle\#}:X\to E$ such that there exists a positive orthomorphism $S\in\Orth(E)$ with
 $|(x,x^{\scriptscriptstyle\#})|\!:=|x^{\scriptscriptstyle\#}(x)|\leq S(\[x\])$ for all $x\in X$; the least $S$
 satisfying the above is denoted by $\[x^{\scriptscriptstyle\#}\]$.

 \Lemma{ 3.4}A linear operator $T:X\to E$ is a member of $X^{\scriptscriptstyle\#}$ if and only if $T$ is norm bounded
 and $\mathbb{B}$-linear; moreover, $\|T\|=\big\|\[T\]\big\|_\infty$, where $\|S\|=\inf\{\lambda>0:\, |S|\leq\lambda I_E\}$.
 \Endproc

 \beginproof~Note first that $\Orth(E)=\mathcal{Z}(E)$ for every normed space $E$ (see \cite[Theorem 3.1.11]{MN}) and
 $\mathcal{Z}(E)$ can be identified with $M$. Consider $T\in L(X,E)$. If $T\in X^{\scriptscriptstyle\#}$,
 $\pi\in\mathbb{P}(E)$, and $x\in X$; then $|T(\pi x)|\leq\pi\[T\]\[x\]$ by Definition 3.3 and so $\pi^\ast T(\pi x)=0$.
 It follows that $\pi Tx=\pi T(\pi x)$. Similarly, $\pi T(\pi^\ast x)=0$ and $T(\pi x)=\pi T(\pi x)$. Thus,
 $T\circ\pi=\pi\circ T$; i.e., $T$ is $\mathbb{B}$-linear. Moreover, $|Tx|\leq\big\|\[T\]\big\|_\infty\[x\]$ so that
 $\|Tx\|=\big\||Tx|\big\|\leq\big\|\[T\]\big\|_\infty\big\|\[x\]\big\|=\big\|\[T\]\big\|_\infty\|x\|$, whence
 $\|T\|\leq\big\|\[T\]\big\|_\infty$. Conversely, if $T$ is norm bounded and $\mathbb{B}$-linear then
 $\|\pi Tx\|_E=\|T(\pi x)\|_E\leq\|T\|\big\|\[\pi x\]\big\|_\infty=\big\|\pi(\|T\|\[x\])|\big\|_E$ for all $\pi\in\mathbb{B}$ and we get
 $|Tx|\leq\|T\|\[x\]$ for all $x\in X$. It follows that $T\in X^{\scriptscriptstyle\#}$ and
 $\|T\|=\big\|\[T\]\big\|_\infty$.~\endproof

 \Lemma{ 3.5}If $\mathcal{X}$ is a Banach space within $\mathbb{V}^{(\mathbb{B})}$ and an order ideal
 $E\subset\mathcal{R}{\downarrow}$ is a Banach lattice then $\mathcal{X}'{\downarrow}\!^{\mathcal{Z}(E)}$
 is $\mathbb{B}$-isometric to $X^{\scriptscriptstyle\#}$. In particular, $X^{\scriptscriptstyle\#}$ is a
 Banach--Kantorovich space over $M$ $($with the obvious identification of $M$ and $\mathcal{Z}(E))$.
 \Endproc

 \beginproof~This follows easily from \cite[Theorem 8.3.6 and Corollary 8.3.7]{DOP} together with
 \cite[Theorem 3.1.11]{MN}.~\endproof

 Assume now that the representation space of  $\mathbb{P}(\mathcal{R}{\downarrow})$ is hyperstone.
 Then there exists an order dense ideal $L\subset\mathcal{R}{\downarrow}$ that is an $AL$-space. Moreover, there is a
 strictly positive order continuous functional $\phi:L\to\mathbb{R}$ such that $\|u\|_L=\phi(|u|)$ for all $u\in L$.
 Put $L^1(\mathbb{B},\phi)\!:=(L,\|\cdot\|_L)$.%

 \Lemma{~3.6}Let $(X,\[\cdot\])$ be a Banach--Kantorovich space over $L^1(\mathbb{B},\phi)$ and $\|\cdot\|\!:=\phi\circ\[\cdot\]$. Then
 $h:x^{\scriptscriptstyle\#}\mapsto\phi\circ x^{\scriptscriptstyle\#}$ is a $\mathbb{B}$-isometry from
 $(X,\[\cdot\])^{\scriptscriptstyle\#}$ onto $(X,\|\cdot\|)'$. In particular, $X'$ is a Banach--Kantorovich space over $M$
 $($with the obvious identification of $M$ and $L^1(\mathbb{B},\phi)')$.
 \Endproc

 \beginproof~Observe first that $X'$ is a Banach--Kantorovich space over $M$ and $\langle x,x'\rangle\leq\phi(\[x\]\[x'\])$
 for all $x\in X$ and $x'\in X'$ (see \cite[Theorem 7.1.4]{DOP}). It is immediate from Lemma 3.4 that $h$ is 
 a~linear
 operator from  $X^{\scriptscriptstyle\#}$ into $X'$. If $\phi\circ x^{\scriptscriptstyle\#}=0$ for some
 $x^{\scriptscriptstyle\#}\in X^{\scriptscriptstyle\#}$, then $\phi(\pi(x,x^{\scriptscriptstyle\#}))=0$ for all $x\in X$ and
 $\pi\in\mathbb{B}$. This implies that $x^{\scriptscriptstyle\#}=0$ as $\phi$ is strictly positive so that $h$ is injective.
 Denote by $B^{\scriptscriptstyle\#}$ and $B'$ the unit balls of $X^{\scriptscriptstyle\#}$ and $X'$, respectively. Then
 $B^{\scriptscriptstyle\#}=\{x^{\scriptscriptstyle\#}\in X^{\scriptscriptstyle\#}:\,\[x^{\scriptscriptstyle\#}\]
 \leq\mathbf{1}\}$ and $B'=\{x'\in X':\,(\forall x\in X)\langle x,x'\rangle\leq\phi(\[x\])\}$. If
 $x^{\scriptscriptstyle\#}\in B^{\scriptscriptstyle\#}$ then $|\langle x,\phi\circ x^{\scriptscriptstyle\#}\rangle|=
 |\phi(\langle x,x^{\scriptscriptstyle\#}\rangle)|\leq\phi(\[x\]\[x^{\scriptscriptstyle\#}\])\leq\phi(\[x\])=\|x\|$ and this
 implies that $h(B^{\scriptscriptstyle\#})\subset B'$. To prove the converse inclusion take $x'\in X'$; i.e.,
 $|\langle x,x'\rangle|\leq\|x'\|\phi(\[x\])$ for all $x\in X$. Since $\phi$ is positive and order continuous, there exists
 a~linear operator $T:X\to L$ such that $x'=\phi\circ T$ and $Tx\leq\|x'\|\[x\]$ for all $x\in X$; see, for example,
 \cite[4.5.2]{SBD}. So $T\in X^{\scriptscriptstyle\#}$ and $h(B')\subset B^{\scriptscriptstyle\#}$. Thus, $h$
 is a $\mathbb{B}$-isometry of $X^{\scriptscriptstyle\#}$ onto $X'$.~\endproof

 \Lemma{~3.7}For every positive order continuous measure $\mu$ on $\mathbb{B}=\mathbb{P}_L(X)$ there exists a unique 
 $N(\mu)\in L^1(\mathbb{B},\phi)$ such that $\mu(b)=\phi(bN(\mu))$ for all $b\in\mathbb{B}$.
 \Endproc

 \beginproof~Identify $\mathbb{B}$ with the Boolean algebra of the components of $\mathbf{1}\in M$, so that
 $\mathbb{B}\subset M$.  Each order continuous measure $\mu:\mathbb{B}\to\mathbb{R}$ admits a unique extension to an order
 continuous functional $f_\mu$ on $M$. Since the order continuous dual $M'_n$ is lattice isometric to $L^1(\mathbb{B},\phi)$, there exists
 exactly one $N(\mu)\in L^1(\mathbb{B},\phi)$ such that $f_\mu(u)=\phi(uN(\mu))$ for all $u\in M$; see \cite[Theorem 3.4.8]{DOP}.~\endproof

   \smallskip

 {\bf Theorem~3.8.}~{\sl Let $X$ be a Banach space with the dual $X'$ and the duality pairing $\langle\cdot,\cdot\rangle$.
 If\/ $\mathbb{B}\!:=\mathbb{P}_L(X)$ then there exists a Banach space $\mathcal{X}$ unique up to linear isometry within
 $\mathbb{V}^{(\mathbb{B})}$ such that the following hold:

 $(1)$~$\mathcal{X}$ has no nontrivial $L$-projections and
 $\mathcal{X}'$ has no nontrivial $M$-projections.

 $(2)$~$X$ is linearly $\mathbb{B}$-isometric to $\mathcal{X}{\downarrow}\!^L$ and $X'$ is linearly $\mathbb{B}$-isometric to
 $\mathcal{X}'{\downarrow}\!^M$.

 $(3)$~There exists a bilinear operator $\langle\!\langle\cdot,\cdot\rangle\!\rangle:X\times X'\to L^1(\mathbb{B},\phi)$ satisfying}
 $$
 \langle\pi x,x'\rangle=\phi(\pi\langle\!\langle x,x'\rangle\!\rangle)\quad(x\in X,\ x'\in X',\ \pi\in\mathbb{B}).
 $$
 \Endproc

 \beginproof~Let $X$ be a Banach space and $\mathbb{B}\!:=\mathbb{P}_L(X)$. As was mentioned after Theorem 2.10, $\mathbb{B}$
 is Bad\'e complete, and so  the function $\mu_x:\mathbb{B}\to\mathbb{R}$ defined as
 $\mu_x(b)\!:=\|bx\|$ $(b\in\mathbb{B})$ is an order continuous measure for every $x\in X$. By Lemma 3.7 there exists a unique $N(\mu_x)$ in $L^1(\mathbb{B},\phi)_+$
 such that $\|bx\|=\phi(bN(\mu_x))$ for all $b\in\mathbb{B}$. Put $\[x\]\!:=N(\mu_x)$ and observe that $\[\cdot\]:X\to L^1(\mathbb{B},\phi)$ is
 a decomposable norm. Indeed, the $d$-decomposability of $\[\cdot\]$ is trivial and $(X,\[\cdot\])$ is $L^1(\mathbb{B},\phi)$-uniformly complete by
 Proposition 2.6. Lemma 2.11 provides the $\mathbb{B}$-completeness of $X$, so that $X$ is a Banach--Kantorovich space by
 Proposition 2.8 and $\|x\|=\big\|\[ x\]\big\|_L$ for all $x\in X$ by definition. So Theorem 3.2 is applicable and guarantees
 the existence within $\mathbb{V}^{(\mathbb{B})}$ of a Banach space $\mathcal{X}$  unique up to linear isometry  whose
 $L^1(\mathbb{B},\phi)$-descent $\mathcal{X}{\downarrow}\!^L$ is $\mathbb{B}$-linearly isometric to $X$. By Lemmas 3.5.and 3.6 $X'$ is linearly
 $\mathbb{B}$-isometric to $\mathcal{X}'{\downarrow}\!^M$ so that 3.8\,(2) holds.

 To prove 3.8\,(1) consider an $L$-projection $\rho\in\mathbb{P}_L(\mathcal{X})$ within $\mathbb{V}^{(\mathbb{B})}$ and
 denote by $P$ the restriction of the descent $\rho{\downarrow}$ to $X$. Using the fact that the descent of the composite
 of mappings within $\mathbb{V}^{(\mathbb{B})}$ is the composition of their descents, we see that $P$ is a projection
 and $\[x\]=\[Px\]+\[(I_X-P)x\]$ for all $x\in X$. Using $\phi$, we get  $\|x\|=\|Px\|+\|(I_X-P)x\|$ $(x\in X)$ so
 that $P\in\mathbb{P}_L(X)$. It follows that there is $b\in \mathbb{B}$ such that $P$ is multiplication by $\chi(b)$; see
 Theorem 3.2. At the same time, by \eqref{CHI}, we have $\chi(b)\in\{0,1\}\subset\mathcal{R}$; see \cite[2.2.6]{KKTop}. It
 follows that $\rho$ is multiplication by $\chi(b)$ and so $\rho$ is trivial. Thus,
 $\mathbb{P}_L(\mathcal{X})=\{0,I_\mathcal{X}\}$ and, by Theorem 2.10, we also have
 $\mathbb{P}_M(\mathcal{X}')=\{0,I_\mathcal{X}\}$.

 Finally, consider the duality $(\mathcal{X},\mathcal{X}')$ and denote by $\delta$ the natural pairing $\delta(x,x')\!:=x'(x)$.
 Then $\delta$ is a bilinear form on $\mathcal{X}\times\mathcal{X}'$ within $\mathbb{V}^{(\mathbb{B})}$ and its
 descent $\delta{\downarrow}$ is a bilinear operator from $\mathcal{X}{\downarrow}\times\mathcal{X}'{\downarrow}$ to
 $\mathcal{R}{\downarrow}$. Using Lemmas 3.4 and 3.5, define a bilinear operator $\langle\!\langle\cdot,\cdot\rangle\!\rangle$
 from $X\times X'$ to $L^1(\mathbb{B},\phi)$ by letting
 $\langle\!\langle x,x'\rangle\!\rangle\!:=\delta{\downarrow}(x,x^{\scriptscriptstyle\#})$ where $x\in X$,
 $x^{\scriptscriptstyle\#}\in X^{\scriptscriptstyle\#}$, and $x'=h(x^{\scriptscriptstyle\#})$. Using $\phi$, we deduce that
 $\phi(\langle\!\langle x,x'\rangle\!\rangle)=\phi(\langle\!\langle x,x^{\scriptscriptstyle\#}\rangle\!\rangle)=
 \langle x,\phi\circ x^{\scriptscriptstyle\#}\rangle=\langle x,x'\rangle$. The proof is complete.~\endproof
 \smallskip

 {\bf Corollary~3.9.}~{\sl Let $X$ be a Banach space with the dual $X'$ and the duality pairing $\langle\cdot,\cdot\rangle$
 and let $\mathbb{B}\!:=\mathbb{P}_L(X)$. Then the following hold:

 $(1)$~$X$ is a Banach--Kantorovich space with mixed norm over $L$ and $\mathcal{P}(X)=\mathbb{P}_L(X)$.

 $(2)$~$X'$ is a Banach--Kantorovich space with mixed norm over $M$ and $\mathcal{P}(X')=\mathbb{P}_M(X')$.

 $(3)$~There exists a bilinear operator $\langle\!\langle\cdot,\cdot\rangle\!\rangle:X\times X'\to L^1(\mathbb{B},\phi)$ satisfying}
 $$
 \langle\pi x,x'\rangle=\phi(\pi\langle\!\langle x,x'\rangle\!\rangle)\quad |\langle\!\langle x,x'\rangle\!\rangle|\leq\[x\]\[x'\],
 \quad(x\in X,\ x'\in X',\ \pi\in\mathbb{B}).
 $$

 \beginproof~This is immediate from Theorems 3.2 and 3.8.~\endproof

 \Remark{~3.10}There is an extensive literature on the idea of continuous (measurable) decomposition (resolution) of the objects of
 functional analysis which stems from the John von Neumann reduction theory; see \cite{Gut, HWW, HK}. One of
 the basic concepts is a \textit{Banach bundle\/} over a topological space called a \textit{base space}. Boolean valued
 analysis is applicable to the same objects whenever the base space is an extremally disconnected compact space. At the
 same time, the principal advantage of this approach is the \textit{transfer principle}: If the object under
 study $X$ admits a Boolean valued representation $\mathcal{X}$, then there is a way to translate ${\rm ZFC}$ theorems on $\mathcal{X}$ to $X$.

 \Remark{~3.11}One of the main reasons for introducing Boolean valued models in functional analysis is to gain insight in
 the structure of a Banach space (or Banach lattice, Banach algebra, etc.) under study by Boolean valued representation.
 Theorem 3.8 is applicable to a wide range of problems in the geometry of Banach spaces whose structure is due
 to the presence of some Boolean algebra of projections. Below we give only two illustrative examples.

 \section{4. Boolean Simplex Spaces}

 In this section we present some conditions on an ordered Banach space that are necessary and sufficient for the space to be a
 predual of an injective Banach lattice. We confine exposition to real Banach spaces, although the method we use works for
 complex spaces as well. The positive cone of an ordered Banach space is assumed closed and the dual space is
 endowed with the dual order.

 \Definition{~4.1}~Let $X$ be a normed vector space ordered by a positive cone $X_+$. Then $X_+$ is said to be
 \textit{normal} if $x\leq y\leq z$ implies that $\|y\|\leq\max\{\|x\|,\|z\|\}$. We say that $X$ \textit{directed\/}
 if the closed unit ball of $X$ is upward directed; i.e. for all $x_1,x_2\in X$ with $\|x_1\|\leq1$ and $\|x_1\|\leq1$
 there exists $y\in X$ such that $x_1,x_2\leq y$ and $\|y\|\leq1$.

 The following characterization of $AL$-spaces is due to Davies \cite{Dav} (also see Asimow \cite{Asi, AE}).

 \Theorem{~4.2}Let $(X, X_+)$ be an ordered Banach space. Then the dual $X'$ is an $AL$-space
 if and only if $X$ has the Riesz decomposition property, $X$ is directed, and $X_+$ is normal.%
 \Endproc

 \Remark{~4.3}~For $L^1$-predual spaces the name \textit{simplex spaces\/} was introduced by Effros \cite{Ef}. Thus, the
 Davies theorem tells us that an ordered Banach space $X$ is a simplex space if and only if $X$ has the Riesz decomposition
 property, $X$ is directed, and $X_+$ is normal. Interpreting this fact in an appropriate Boolean valued model we will
 characterize the preduals of injective Banach lattices. It is natural to call them \textit{Boolean simplex spaces}.
 The following result on Boolean valued representation of injective Banach spaces together with Theorem 3.8 provides a
 key to the characterization of Boolean simplex spaces.
 \smallskip

 Given a Banach lattice $X$, denote by $\mathbb{M}(X)$ the set of the $M$-projections that are simultaneously band
 projections: $\mathbb{M}(X)\!:=\mathbb{P}(X)\cap\mathbb{P}_M(X)$. Then $\mathbb{M}(X)$ is a subalgebra of the Boolean
 algebra $\mathbb{P}(X)$. Moreover, $\mathbb{M}(X)$ is an order closed subalgebra of $\mathbb{P}(X)$ whenever each upward
 directed set in the unit ball of $X$ has the least upper bound belonging to the unit ball of $X$; see \cite{Hay}.%

 \Theorem{~4.4}A Banach lattice $X$ is injective with $\mathbb{B}=\mathbb{M}(X)$ if and only if there exists an $AL$-space
 $\mathcal{X}$ unique up to lattice isometry within $\mathbb{V}^{(\mathbb{B})}$ whose bounded descent
 $\mathcal{X}{\downarrow}\!^M$ is $\mathbb{B}$-isometric to $X$.
 \Endproc

 \beginproof~See \cite[Theorem 2.1]{Kus_IBL} or \cite[Theorem 5.9.1]{KKTop}.~\endproof

 \Lemma{ 4.5}$\mathbb{P}_M(X)=\mathbb{M}(X)$ for every  injective Banach lattice $X$.
 \Endproc

 \beginproof~Let $X$ and $\mathcal{X}$ be the same as in Theorem 4.4. The mapping
 $\pi\mapsto\pi{\Downarrow}\!:=\pi{\downarrow}|_X$ is an isomorphism of $\mathbb{P}_M(\mathcal{X}){\downarrow}$
 and $\mathbb{P}_M(X)$. This can be proved as in \cite[Theorem 5.9.1]{KKTop}. If $P\in\mathbb{P}_M(X)$ and
 $P\notin\mathbb{M}(X)$ then $\pi\!:=P{\uparrow}$ is a nontrivial $M$-projection in $\mathcal{X}$ according to
 \cite[Theorem 5.8.12]{KKTop}, because $\mathbb{B}$-linearity of $P$ amounts to saying that $P$ commutes with each operator
 from $\mathbb{M}(X)$. Apart from the trivial case $X=\mathbb{R}$, this contradicts to Behrends' dichotomy
 (see \cite[Theorem 1.8]{HWW}), since the $AL$-space $\mathcal{X}$ has nontrivial $L$-projections. As to the trivial case,
 the $M$-projections of $\mathcal{R}{\downarrow}\!^M\simeq C(K)$ coincide with the characteristic projections
 $P_U(f) = \chi_Uf$ $(f\in C(K))$ for clopen sets $U\subset K$ \cite[Example 1.4\,(a)]{HWW}.~\endproof

 Now we introduce the Boolean versions of the concepts of normality and directedness of ordered Banach spaces which appear
 in the Davies characterization of $L^1$-predual spaces.

 \Definition{~4.6}~The norm on $X$ is said to be $\mathbb{B}$-\textit{normal} if for $x,y,z\in X$ with $x\leq y\leq z$ there
 exists a projection $\pi_0\in\mathbb{B}$ such that $\|\pi_1y\|\leq\|\pi_1x\|$ and $\|\pi_2 y\|\leq\|\pi_2 z\|$ for all
 $\pi_1,\pi_2\in\mathbb{B}$ with $\pi_1\leq\pi_0$ and $\pi_2\leq\pi_0^\ast$.

 \Definition{~4.7}~Given $a\in X$, define $B_a(X)\subset X$ as the set of $x\in X$ with $\|\pi x\|\leq\|\pi a\|$
 for all $\pi\in\mathbb{P}_L(X)$. Say that $B_a(X)$ is upward $\mathbb{B}$-\textit{directed\/} if, for all
 $x_1,x_2\in B_a(X)$, there exists $y\in B_a(X)$  with $y\geq x_1,x_2$.

 \Lemma{ 4.8}Let $X$ be the $E$-descent of an ordered Banach space $\mathcal{X}\in\mathbb{V}^{(\mathbb{B})}$ where $E$ is a
 Dedekind complete Banach lattice and $\mathbb{B}=\mathbb{P}(E)$. Then the following   hold:

 $(1)$~$\mathcal{X}_+$ is normal if and only if $X_+$ is $\mathbb{B}$-normal.

 $(2)$~$\mathcal{X}$ is directed if and only if $B_a(X)$ is upward $\mathbb{B}$-directed for all $a\in X$.
 \Endproc

 \beginproof~Recall that the $E$-valued norm $\[\cdot\]$ of $X$ is the descent of the norm $\|\cdot\|_{\mathcal{X}}$ of
 $\mathcal{X}$ and $\|x\|_X=\big\|\[x\]\big\|_\infty$ for all $x\in X$. Observe first that for $a,b\in E_+$ we have
 $a\leq b$ if and only if $\|\pi a\|_E\leq\|\pi a\|_E$ for all $\pi\in\mathbb{B}$ (see the proof of Lemma 2.8). It follows
 that the conditions $\[x\]\leq\[y\]$ and $\|\pi x\|\leq\|\pi y\|$ for all $\pi\in\mathbb{B}$
 are equivalent. Now it is clear that $B_a(X)=\{x\in X:\ \[x\]\leq\[a\]\}$.

 $(1)$~Suppose that $\mathcal{X}_+$ is normal and pick $x,y,z\in X$ with $x\leq y\leq z$. Then
 $[\![x\leq y\leq z]\!]=\mathbbm{1}$ and $\|y\|_\mathcal{X}\leq\max\{\|x\|_\mathcal{X},\|z\|_\mathcal{X}\}$. It follows
 that $\[y\]\leq\[x\]\vee\[z\]$ and there exists a band projection $\pi_0\in\mathbb{B}$ such that
 $\[x\]\vee\[z\]=\pi_0\[x\]+\pi_0^\ast\[z\]$. The simple argument indicated above shows that
 $\|\pi y\|\leq \|\pi\pi_0 x\|+\|\pi\pi_0^\ast z\|$ for all $\pi\in\mathbb{B}$. Taking $\pi\!:=\pi_1\leq\pi_0$, we get
 $\|\pi_1y\|\leq\|\pi_1x\|$ and for $\pi\!:=\pi_2\leq\pi_0^\ast$ we have $\|\pi_2 y\|\leq\|\pi_2 z\|$. The converse can
 be proven by similar arguments.

 $(2)$~Denote by $B_1(\mathcal{X})$ the closed unit ball of $\mathcal{X}$ and put $B_\alpha(\mathcal{X})=\alpha B_1(\mathcal{X})$.
 Note that $B_1(\mathcal{X})$ is upward directed if and only if so is $B_\alpha(\mathcal{X})$ for all
 $\alpha=\|a\|_\mathcal{X}\in\mathcal{R}$ with $a\in\mathcal{X}$.
 The straightforward calculation of Boolean truth-values yields that $B_\alpha(\mathcal{X})$ is upward directed if and only if
 $B_a(X)$ is upward directed for $\alpha=\[a\]$. Thus the directedness of $X$ implies that $B_a(X)$ is upward directed for
 all $a\in X$, since $B_\alpha(\mathcal{X}){\downarrow}\subset X$. To show the converse, choose a~partition of unity
 $(\pi_\xi)$ in $\mathbb{B}$ and a family $(a_\xi)$ in $X$ such that $\pi_\xi a_\xi\in X$ and
 $\pi_\xi\[a_\xi\]=\pi_\xi\mathbf{1}\in X$ for all $\xi$. For $\mathbb{O}\ne b\in\mathbb{B}$ denote by $b\mathbb{B}$
 the relative subalgebra $[\mathbb{O},b]$. Assuming that all $B_{a_\xi}(X)$ are upward directed and taking into account
 the rule for transition to the relative universe $\mathbb{V}^{(b\mathbb{B})}$ (see \cite[1.3.7]{KKTop}) we infer that
 $B_{\alpha_\xi}(\mathcal{X})=B_1(\mathcal{X})$ within $\mathbb{V}^{(\pi_\xi\mathbb{B})}$ and so $B_1(\mathcal{X})$
 is upward directed within $\mathbb{V}^{(\mathbb{B})}$.~\endproof


 \Theorem{~4.9}For an ordered Banach space $X$ the following are equivalent:

 {\rm (1)}~$X'$ is an injective Banach lattice. 

 {\rm (2)}~The four conditions hold: 

 $\quad\quad (a)$~$X$ has the Riesz  decomposition property.

 $\quad\quad (b)$~$\mathbb{P}_L(X)$ consists of positive projections.

 $\quad\quad (c)$~$X_+$ is\/ $\mathbb{P}_L(X)$-normal.

 $\quad\quad (d)$~$B_a(X)$ is upward $\mathbb{P}_L(X)$-directed for all $a\in X$.

 In cases (1) and/or (2) the Boolean algebras $\mathbb{M}(X')$ and $\mathbb{P}_L(X)$ are isomorphic.%
 \Endproc

 \beginproof~Let $\mathcal{X}$ be a Boolean valued representation of $X$ as in Theorem 3.8. Just as in
 \cite[Theorem 5.9.1]{KKTop} we can demonstrate that $X_+{\uparrow}$ is a pointed cone in $\mathcal{X}$
 over $\mathbb{R}_+^{\scriptscriptstyle\wedge}$; i.e. $X_+{\uparrow}+X_+{\uparrow}\subset X_+{\uparrow}$
 and $\mathbb{R}^{\scriptscriptstyle\wedge}_+\cdot X_+{\uparrow}\subset X_+{\uparrow}$. Define the positive cone
 $\mathcal{X}_+$ as the closure of $X_+{\uparrow}$. Then $\mathcal{X}_+$ is a closed, possibly unpointed, cone
 so that $(\mathcal{X},\mathcal{X}_+)$ is a preordered Banach space within $\mathbb{V}^{(\mathbb{B})}$. Recall that
 for $A\subset X$ we have $A{\uparrow\downarrow}=\mix(A)\subset\mathcal{X}{\downarrow}$, where $\mix(A)$ consists of
 all $x\in\mathcal{X}{\downarrow}$ for which there exists a partition of unity $(\pi_\xi)$ in $\mathbb{B}$ and a
 family $(a_\xi)$ in $A$ such that $\pi_\xi x=\pi_\xi a_\xi$ for all $\xi$; see \cite[1.6.6]{KKTop}. Moreover, if
 $x\in X$ then $x=\sum_{\xi\in\Xi}\pi_\xi a_\xi$ by Lemma 2.11.~\smallskip

 $(1)\Longrightarrow(2)$:~If $X'$ is an injective Banach lattice then $\mathbb{P}_M(X')=\mathbb{M}(X')$
 by Lemma 4.4 and $\mathbb{P}_L(X)$ is isomorphic to $\mathbb{M}(X')$ by Theorem 2.10\,(3).
 Applying the obvious duality representation $\mathbb{P}_L(X)=\{\pi\in\mathcal{L}(X):\,\pi'\in\mathbb{P}_M(X)\}$,
 we conclude that $L$-projections on $X$ are positive and so $X_+$ is invariant under all $L$-projections.
 Together with Lemma 2.11, this implies that $X_+=\mix(X_+)\cap X=X_+{\uparrow\downarrow}\cap X$, which in turn yields that
 $X_+{\uparrow}$ is closed in $\mathcal{X}$; i.e. $\mathcal{X}_+=X_+$. Thus, $\mathcal{X}$ is an ordered
 Banach space and $\mathcal{X}'$ is an $AL$-space within $\mathbb{V}^{(\mathbb{B})}$ by Theorem 4.4. The Boolean
 valued transfer principle enables us to apply the Davies Theorem to $\mathcal{X}$ and state that $\mathcal{X}$
 has the Riesz decomposition property, $\mathcal{X}$ is directed, and $\mathcal{X}_+$ is normal. By Lemma 4.7 $B_a(X)$
 is upward $\mathbb{P}_L(X)$-directed for all $a\in X$ and $X_+$ is\/ $\mathbb{P}_L(X)$-normal. The Riesz
 decomposition property for $\mathcal{X}$ can be written in the equivalent form $\Sigma([0,a]\times[0,b])=[0,\Sigma(a,b)]$,
 where $a,b\in\mathcal{X}_+$ and $\Sigma$ is addition on $\mathcal{X}$. Then $\Sigma{\downarrow}$ is addition on
 $X$ and $\Sigma([0,a]\times[0,b]){\downarrow}=\Sigma{\downarrow}([0,a]{\downarrow}\times[0,b]{\downarrow})$. The latter
 is equivalent to the Riesz decomposition property for $X$, since obviously the descent of an internal order interval
 $\{x\in\mathcal{X}:\,0\leq x\leq a\}$ coincides with $\{x\in X:\,0\leq x\leq a\}$.
~\smallskip

 $(2)\Longrightarrow(1)$:~Assume that $X$ satisfies the conditions $(a)$--$(d)$ of (2) and put $\mathbb{B}=\mathbb{P}_L(X)$.
 By $(a)$, $X_+$ is invariant under $L$-projections, which together
 with the closedness of $X_+$ implies that $X_+$ is $\mathbb{B}$-complete. Defining the positive cone in $\mathcal{X}$ by
 $\mathcal{X}_+\!:=X_+{\uparrow}$, we can demonstrate just as in \cite[Theorem 5.9.1]{KKTop} that $(\mathcal{X},\mathcal{X}_+)$
 is an ordered Banach space within $\mathbb{V}^{(\mathbb{B})}$. As was shown above, $X$ has the Riesz decomposition property if
 and only if so is $\mathcal{X}$ within $\mathbb{V}^{(\mathbb{B})}$. Prove that $\mathcal{X}$ is directed and $\mathcal{X}_+$
 is normal. By Lemma 4.7, $(c)$ and $(d)$ imply that $\mathcal{X_+}$ is normal and $\mathcal{X}$ is directed. By the Boolean
 valued transfer principles Theorem 4.2 is true within $\mathbb{V}^{(\mathbb{B})}$ and so $\mathcal{X}$ is an $L^1$-predual space; i.e. $\mathcal{X}'$ is an $AL$-space. It remains to appeal again to Theorem 4.4. By Lemma 4.5
 $\mathbb{M}(X')$ and $\mathbb{P}_L(X)$ are isomorphic.~\endproof

 \Remark{~4.10}As was shown by Asimow in \cite[Theorem 1]{Asi}, directedness of the unit ball is also involved in
 characterization of those ordered Banach spaces whose dual cones are \textit{well-capped}. This concept, dating back to
 Choquet, is of interest due to the results of the following type: A closed convex well-capped subset of a locally convex
 topological vector space is the closed convex hull of its extreme points and extremal rays; see Asimow
 \cite[Theorem 2.2]{Asim}. When studying the extremal structure of sets of operators, the ordinary caps are of little use.
 The theory of  \textit{operator caps} and \textit{well-capped set of
 operators} was built in \cite{Kut3} in the spirit of the Choquet cap theory; see also \cite[Chapter 2, \S~5]{SBD}.

 \section{5. Centerable Subsets}

 Now we characterize preduals of injective Banach lattices in  terms of centerable sets. Let $X$ be a Banach space,
 let $\mathbb{B}$ be a Boolean algebra of projections on $X$, and let $A$ be a bounded subset of $X$.

 \Definition{~5.1}The $\pi$-\textit{diameter} $\delta_\pi(A)$ and the \textit{Chebyshev $\pi$-radius}
 $r_\pi(A)$ of $A$ are defined as $\delta_\pi(A)=\sup\{\|\pi(a- b)\|:\, a, b\in A\}$ and $r_\pi(A)=\inf_{x\in X}r_\pi(A,x)$
 where $r_\pi(A, x)=\sup_{a\in A}\|\pi(x-a)\|$ for all $x\in X$. It is easily seen that $\delta_\pi(A)\leq2r_\pi(A)$.
 If $\pi$ coincides with the identity operator $I_X$, then we write $\delta(A)$ and $r(A)$ instead of $\delta_\pi(A)$
 and $r_\pi(A)$, respectively. If $\delta(A)=2r(A)$, then A is said to be \textit{centerable}.

 The next result is due to Duan and Lin \cite[Theorem 2.7]{DL}.

 \Theorem{~5.2}For a real Banach space $X$ the following are equivalent:

 $(1)$~$X$ is an $L^1$-predual space.

 $(2)$~Every four-point subset of $X$ is centerable.

 $(3)$~Every finite subset of $X$ is centerable.

 $(4)$~Every compact subset of $X$ is centerable.
 \Endproc

 \Remark{~5.3}Theorem 5.2 is true also for a complex Banach spaces as was proved in \cite[Theorem 2.13]{DL}. Moreover,
 it is pointed out in \cite[Remark 2.14]{DL} that this result cannot be sharpened anymore; i.e., the centerability of
 every three-point subset of a real or complex Banach space $X$ does not imply that $X$ is an $L^1$-predual space.

 Denote the set of all partitions (respectively, all countable partitions) of unity in $\mathbb{B}$ by $\Prt(\mathbb{B})$
 (respectively, $\Prt_\sigma(\mathbb{B}))$.

 \Definition{~5.4}~A subset $A\subset X$ of a $\mathbb{B}$-complete Banach space $X$ is called
 $\mathbb{B}$-\textit{centerable\/} if $2r_\pi(\mix(A))=\delta_\pi(\mix(A))$ for all $\pi \in\mathbb{B}$, where
 $$
 \mix(A)\!:=\Big\{\sum\nolimits_{\xi\in\Xi}\pi_\xi a_\xi\in X:\
 (\pi_\xi)_{\xi\in\Xi}\in\Prt\nolimits_\sigma(\mathbb{B}),\
 \{a_\xi:\ \xi\in\Xi\}\subset A\Big\}.
 $$

  \Definition{~5.5}~Suppose that $X$ is a $\mathbb{B}$-cyclic Banach lattice. Say that a~sequence $(x_n)_{n\in\mathbb{N}}$
 in $X$\/ {\it $\mathbb{B}$-approximates~$x\in X$} if, for each $k\in \mathbb{N}$, we have
 $$
 \inf\{\sup_{n\geqslant k}\|\pi_n(x_n-x)\|:\ (\pi_n)_{n\geq
 k}\in\Prt_\sigma(\mathbb{B})\}=0.
 $$ Call a set
 $K\subset X$ {\it $\mathbb{B}$-compact} (or \textit{cyclically compact}, \cite[2.12.C]{KKTop}) if $K$ is
 $\mix$-complete and for every sequence  $(x_n)_{n\in\mathbb{N}}\subset K$ there is $x\in K$ such that
 $(x_n)_{n\in\mathbb{N}}$  $\mathbb{B}$-approximates~$x$.

 We are now ready to state the main result of this section.

 \Theorem{~5.6}For a real Banach space $X$ the following are equivalent:

 $(1)$~$X'$ is an injective Banach space with $\mathbb{M}(X')$ isomorphic to\/ $\mathbb{B}\!:=\mathbb{P}_L(X)$.

 $(2)$~Every four-point subset of $X$ is\/ $\mathbb{B}$-centerable.

 $(3)$~Every finite subset of $X$ is\/ $\mathbb{B}$-centerable.

 $(4)$~Every $\mathbb{B}$-bounded $\mix$-compact subset of $X$ is\/ $\mathbb{B}$-centerable.

 $(5)$~For every $\mix$-compact subset $A$ of $X$ there exists a partition of unity $(\pi_\xi)_{\xi\in\Xi}$ in
 $\mathbb{B}$ such that $\pi_\xi A$ is\/ $\pi_\xi\mathbb{B}$-centerable in $\pi_\xi X$for all $\xi\in\Xi$.
 \Endproc
 \smallskip

 The technical details of the proof are highlighted in the following two lemmas:

 \Lemma{~5.7}Let $\mathcal{X}$ be a Banach space within $\mathbb{V}^{(\mathbb{B})}$. Put $X\!=\mathcal{X}{\downarrow}$
 and consider a subset $A\subset X$ with $\mix(A)=A$ and $\[A\]$ order bounded in $\mathcal{R}{\downarrow}$. Then
 $[\![\delta(A{\uparrow})=2r(A{\uparrow})]\!]=\mathbbm{1}$ if and only if $2R(A)=\Delta(A)$ in $\mathcal{R}{\downarrow}$,
 where
 \begin{gather*}
 \Delta(A)\!:=\sup\{\[a- b\]:\, a, b\in A\},
 \\
 R(A)\!:=\inf\nolimits_{x\in X}R(x,A),\quad R(x,A)\!:=\sup\nolimits_{a\in A}\[x-a\].
 \end{gather*}
 \Endproc

 \beginproof~The order boundedness of $\[A\]$ guarantees the existence of $\Delta(A)$ and $R(A)$ in
 $\mathcal{R}{\downarrow}$. The formula $\delta(A{\uparrow})=2r(A{\uparrow})$ is equivalent to $2r(A{\uparrow})\leq\delta(A{\uparrow})$,
 while the latter can be rewritten as
 \begin{equation}\label{Eq}
 \Psi\equiv(\forall n\in\mathbb{N}^{\scriptscriptstyle\wedge})(\exists x\in\mathcal{X})(\exists a,b\in A{\uparrow}) (\forall c\in A{\uparrow})\ 2\|x-c\|_\mathcal{X}\leq\|a-b\|_\mathcal{X}+1/n.
 \end{equation}
 Calculating the Boolean truth-values for the quantifiers in the formula $[\![\Psi]\!]=\mathbbm{1}$ (using the rules
 \cite[1.4.5\,(1), 1.6.2, and 1.6.6]{KKTop}), we arrive at the equivalent  condition
 $$
 (\forall n\in\mathbb{N})\ (\exists x\in\mathcal{X}{\downarrow})\ (\exists a,b\in A)\
 (\forall c\in A)\ 2\[x-c\]\leq\[a-b\]+(1/n)\mathbf{1}.
 $$
 It follows that $\sup_{c\in A}2\[x-c\]\leq\[a-b\]+(1/n)1^{\scriptscriptstyle\wedge}$ and so
 $\sup_{c\in A}2\[x-c\]\leq\Delta(A)+(1/n)1^{\scriptscriptstyle\wedge}$,
 whence we get the inequality $2R(A)\leq\Delta(A)$ equivalent to $2R(A)=\Delta(A)$. Conversely, if the latter holds then
 for every $n\in\mathbb{N}$ there exist a partition of unity $(\pi_\xi)_{\xi\in\Xi}$ in $\mathbb{B}$ and families
 $(a_\xi)_{\xi\in\Xi}$ and $(b_\xi)_{\xi\in\Xi}$ in $A$ and $(x_\xi)_{\xi\in\Xi}$ in $X$ such that
 $2\pi_\xi\[x_\xi-c\]\leq\[a_\xi-b_\xi\]+(1/n)\mathbf{1}$ for all $c\in A$. Put
 $x=\mix_{\xi\in\Xi}\pi_\xi x_\xi$, $a=\mix_{\xi\in\Xi}\pi_\xi a_\xi$, and $b=\mix_{\xi\in\Xi}\pi_\xi b_\xi$
 and observe that $a,b\in A$, $x\in\mathcal{X}{\downarrow}$ and  $2\pi_\xi\[x-c\]\leq\[a-b\]+(1/n)\mathbf{1}$
 for all $c\in A$. The latter is equivalent to \eqref{Eq}, as can be checked by the direct calculation
 of Boolean truth-values.~\endproof

 \Lemma{~5.8}If $\pi\in\mathbb{B}$ and $x\in X$ then $\{\[\pi(a-b)\]:\, a, b\in\pi A\}$ and $\{\[\pi(x-c)\]:\, c\in\pi A\}$
 are upward directed, whereas $\{R(\pi x,\pi A):\, x\in X\}$ is
 downward directed.
 \Endproc

 \beginproof~Since $\pi R(x, A)=R(\pi x, \pi A)$ and $\pi D(A)=D(\pi A)$, we can assume without loss
 of generality that $\pi=I_X$. Since the first two claims can be checked similarly, we consider only the first one.
 If $a,b\in A$ and $a',b'\in A$ then we can pick a band projection $\pi_0\in\mathbb{B}$ such that
 \begin{equation*}
 \[a-b\]\vee\[a'-b'\]=\pi_0\[a-b\]+\pi_0^\ast\[a'-b'\]=\[(\pi a+\pi_0^\ast a')-(\pi b+\pi_0^\ast b')\]=\[u-v\]
 \end{equation*}
 where  $u\!:=\pi a+\pi_0^\ast a'\in A$ and $v\!:=\pi b+\pi_0^\ast b'\in  A$ and the second equality is due to
 Proposition 2.3. Next, take  $x,y\in X$ and choose $\pi_0\in\mathbb{B}$ so that $R(x,A)\wedge
 R(y,A)=\pi_0R(x,A)+\pi_0^\ast R(y,A)$. Then
 \begin{align*}
 R(x,A)\wedge R(y,A)&=R(\pi_0 x,\pi_0 A)+R(\pi_0^\ast y,\pi_0^\ast  A)
 \\
 &=\sup_{a\in A}\[\pi_0(x-a)\]+\sup_{a\in A}\[\pi_0^\ast(y-a)\]
 \\
 &\geq\sup_{a\in A} \[(\pi_0 x+\pi_0^\ast y)-a\]
 \\
 &=R(u,A);
 \end{align*}
 here $u\!=\pi_0 x+\pi_0^\ast y\in A$.~\endproof
 \smallskip

 {\sc Proof of Theorem 5.6.}~The implications $(5)\Longrightarrow(4)\Longrightarrow(3)\Longrightarrow(2)$ are obvious.
 The rest of the proof is based on Theorem 3.8, according to which there is no loss of generality in assuming that
 $X=\mathcal{X}{\downarrow}\!^L$ and $X'=\mathcal{X'}{\downarrow}\!^M$ for some Banach space $\mathcal{X}$ within
 $\mathbb{V}^{(\mathbb{B})}$, where $\mathbb{B}=\mathbb{P}_L(X)$.

 $(1)\Longrightarrow(5)$:~If $\mathcal{X}'$ is an $AL$-space then by the Boolean valued transfer principle and Theorem 5.2
 every compact subset of $\mathcal{X}$ is centerable. Take a cyclically compact subset $A\subset X$. By
 \cite[2.14.C.3]{KKTop} (or \cite[8.5.1]{DOP}) $A{\uparrow}$ is a compact subset of $\mathcal{X}$ and so
 $[\![\delta(A{\uparrow})=2r(A{\uparrow})]\!]=\mathbbm{1}$. Moreover, $[\![A{\uparrow}$ is norm bounded$]\!]=\mathbbm{1}$
 and therefore $\[A\]$ is order bounded in $\mathcal{R}{\downarrow}$. By Lemma 5.7 $2R(A)=\Delta(A)$ but it may
 happen that $2R(A)$ and $\Delta(A)$ are not in $L^1(\mathbb{B},\phi)$. Choose a partition of unity $(\pi_\xi)_{\xi\in\Xi}$
 such that $2\pi_\xi R(A),\pi_\xi\Delta(A)\in L^1(\mathbb{B},\phi)$ for all $\xi\in\Xi$. Thus for every $\xi\in\Xi$
 we have the equation $2\pi_\xi R(A)=\pi_\xi\Delta(A)$ which in turn, by the argument of the proof of Lemma 2.8,
 amounts to stating that $\|2\pi\pi_\xi R(A)\|_L=\|\pi\pi_\xi\Delta(A)\|_L$ for all $\pi\in\mathbb{B}$. Considering the
 order continuity of the $AL$-space $L^1(\mathbb{B},\phi)$ and using Lemma 5.8, for every $\pi\leq\pi_\xi$ we deduce that
 \begin{multline}\label{OC}
 2r_\pi(A)=2\inf_{x\in X}\sup_{c\in A}\|\pi(x-c)\|=2\Big\|\inf_{x\in X}\sup_{c\in A}\,\[\pi(x-c)\]\,\Big\|_L
 =\|2\pi R(A)\|_L
 \\
 =\|\pi\Delta(A)\|_L=\Big\|\sup_{a,b\in A}\[\pi(a-b)\]\Big\|_L=\sup_{a,b\in A}\|\pi(a-b)\|=\delta_\pi(A).
 \end{multline}

 $(2)\Longrightarrow(1)$:~Assume  that every four-point subset of $X$ is $\mathbb{B}$-centerable and let $A$ be
 a~four-point
 subset of $\mathcal{X}$; i.e.,
 \begin{multline*}
 [\![(\exists x_1,x_2,x_3,x_4\in\mathcal{X})(\forall x\in\mathcal{X})\big(x\in A
 \\
 \leftrightarrow (x_1\in A)\vee(x_2\in A)\vee(x_3\in A)\vee(x_4\in A)\big)]\!]=\mathbbm{1}.
 \end{multline*}
 Using the maximum principle, we conclude that there are
 $x_1,x_2,x_3,x_4\in\mathcal{X}{\downarrow}$ satisfying $A{\downarrow}=\mix(\{x_1,x_2,x_3,x_4\})$. There exists a
 partition of unity $(\pi_\xi)_{\xi\in\Xi}$ in $\mathbb{B}$ such that $\pi_\xi x_j\in X$ for all $\xi\in\Xi$ and
 $j=1,2,3,4$ and so $\pi_\xi A{\downarrow}\subset X$ for all $\xi\in\Xi$. By hypothesis,
 $2r_{\pi\pi_\xi}(A{\downarrow})=\delta_{\pi\pi_\xi}(A{\downarrow})$ for all $\xi\in\Xi$
 and $\pi\in\mathbb{B}$. Arguing as in the proof of \eqref{OC}, we see that
 $2\pi_\xi R(A{\downarrow})=\pi_\xi\Delta(A{\downarrow})$ for all $\xi\in\Xi$; and so
 $2R(A{\downarrow})=\Delta(A{\downarrow})$. Now, working within $\mathbb{V}^{(\mathbb{B})}$ and using
 Lemma 5.7 and formula $A=A{\downarrow\uparrow}$, we infer that $[\![\delta(A{\uparrow})=2r(A{\uparrow})]\!]=\mathbbm{1}$.
 Thus, according to the transfer principle, every four-point subset of $\mathcal{X}$ is\/ centerable and $\mathcal{X}'$
 is an $AL$-space by Theorem 5.2. It remains to appeal to Theorem 4.4.~\endproof

 \Remark{~5.9}Many interesting geometric characterizations of Banach spaces are due to various ball intersection properties; see
 \cite{HWW, Lin}. In \cite[Theorem 6.1]{Lin} Lindenstrauss proved that the dual $X'$ of a real Banach space $X$ is isometric
 to an $L^1$-space if and only if $X$ has the $(4,2)$-intersection property; i.e., every collection of four mutually
 intersecting closed balls has  nonempty intersection. To formulate the Boolean version of this result, we define a
 $\mathbb{B}$-\textit{cell\/} as a set of the form
 $B(a,r)\!:=\{x\in X:\,\|\pi(x-a)\|\leq\|\pi r\|_L\text{~for all~}\pi\in\mathbb{B}\}$, where $\mathbb{B}\!:=\mathbb{P}_L(X)$,
 $a\in X$, and $r\in L^1(\mathbb{B},\phi)$. The above technique allows us to state that the dual $X'$
 of a real Banach space $X$ is an injective Banach lattice with $\mathbb{M}(X')$ isomorphic to $\mathbb{B}\!:=\mathbb{P}_L(X)$
 if and only if every collection of four mutually intersecting $\mathbb{B}$-cells in $X$ has  nonempty intersection.

 \section{6. Concluding Remarks}

 {\bf 6.1.}~Theorem 4.4 yields the following representation result: Given an injective Banach lattice $X$ with
 $\mathbb{B}=\mathbb{M}(X)$, there exists a strictly positive Maharam operator $\Phi:X\to M$ (i.e., an interval
 preserving order continuous linear operator) with the Levi property, where $M$ is a Dedekind complete $AM$-space
 with unit and $\mathbb{P}(M)$ isomorphic to $\mathbb{B}$ such that $\|x\|=\|\Phi(|x|)\|_\infty$ for all $x\in X$; see
 \cite[Corollary 4.5]{Kus_IBL}. The concept of \textit{Maharam operator\/} stems from the articles \cite{Mah1, Mah2}
 on the representation of positive operators on function spaces; the term was introduced in \cite{LS}. Luxemburg was
 the first to appreciate Maharam's contribution. In his joint articles with Schep~\cite{LS} and de~Pagter \cite{LuP}
 some~portion of Maharam's theory was extended to positive operators in Dedekind complete vector lattices.
 \smallskip

 {\bf 6.2.}~Luxemburg was a pioneer and promoter of blending model theory and functional analysis. He pointed out that the
 Maharam operators may play a fundamental role not only in the theory of positive operators but also in Boolean valued
 analysis. His article \cite{Lux} in the Maharam anniversary volume states:%
 \medskip

 \textit{Finally we like to mention that the $F$-measure algebras introduced by Maharam in \cite{Mah1} have also recently
 appeared in the literature in the form of Boolean-valued models of standard numerical measure algebras. It may be of
 interest to explore further the properties of such abstract $F$-measure algebras from the point of view of the theory
 developed by D.~Maharam. For further details we have to refer the reader to {\rm \cite{Tak}}.}
  \medskip

 The development of Maharam's ideas within Boolean valued analysis is due to other authors. Thanks to the articles
 \cite{Lux, LuP, LS}, these ideas were fruitfully implanted into operator theory,
 convex analysis, and elsewhere; see   \cite{DOP, SBD, BVA, KKTop}.
 \smallskip


 \Address{\textit{Received October ??, 2019.}
 \\[6pt]
 \textsc{Kusraev A.G.}\\
 \textsc{Southern Mathematical Institute}\\
 \textsc{Vladikavkaz Science Center of the Russian Academy of Sciences} \\
 \textsc{Khetagurov North Ossetian State University, Vladikavkaz, Russia}\\
 \textit{E-mail address:}  \verb"kusraev@smath.ru"
 \\[4pt]
 \textsc{Kutateladze S. S.} \\
 \textsc{Sobolev Institute of Mathematics}\\
 \textsc{Siberian Branch of the Russian Academy of Sciences, Novosibirsk, Russia}\\
 \textit{E-mail address:} \verb"sskut@math.nsc.ru"}

 \end{document}